# Prime Identification and Composite Filtering Using GM-(n+1) Sequences

**Author:**
Dr. Fadwa Hamdi Barakat
Jadara University, Department of Design Communication,
Faculty of Languages & Literature, Irbid B.O.Box 733,
Irbid, 21110, Jordan
Email: fadwa@jadara.edu.jo

Abstract

This paper presents a distinctive prime detection approach. This method use GM-(n+1) sequences to effectively eliminate complex numbers. The sequences, which consist of odd a number of (n+1), exclude all components except for the initial prime integer. Only the first prime number is presented. This research proposes an approach using this model to identify exceptional candidates and examine their distribution. This study examines the interconnections among the laws of division, basic gaps, and their applications in analytical procedures. Computer studies may provide a novel perspective on the theory of prime numbers, demonstrating the effectiveness of this approach in refining the search space for primes.

## 1. Introduction
All topics of number theory, cryptography, numerical analysis, and algorithm creation consider prime numbers to be vital. Candall and Pomerance (2005) and Agrawal, Kayal, and Saxena (2004) cite the experimental section, the sieve of Eratosthenes, and probabilistic tests such as the Miller-Rabin as not only effective but also inexpensive. Furthermore, the filter composites method provides a composite sifting method that takes advantage of the arithmetic sequence structure in the GM-(n+1) sequence approach. This study attempts to explore the design, computational execution, and the theoretical impacts of this new approach regarding its effectiveness compared to other existing ones (Hardy & Wright, 2008; Tenenbaum, 1995).

## 2. Definition of GM-(n+1) Sequences
For a positive integer, n, the complete set of GM-(n+1) sequences is exhaustive.

The graph of GM-(n+1) is as given below: $\{(n+1)k \mid k \in \mathbb{Z}^+\}$.

Take GM-3 for example: {3, 6, 9, 12, 15, … } The only prime number is 3.

- GM-5: {5, 10, 15, 20, 25, … } Only 5 is prime.

- GM-7: {7, 14, 21, 28, 35, … } Only 7 is prime.

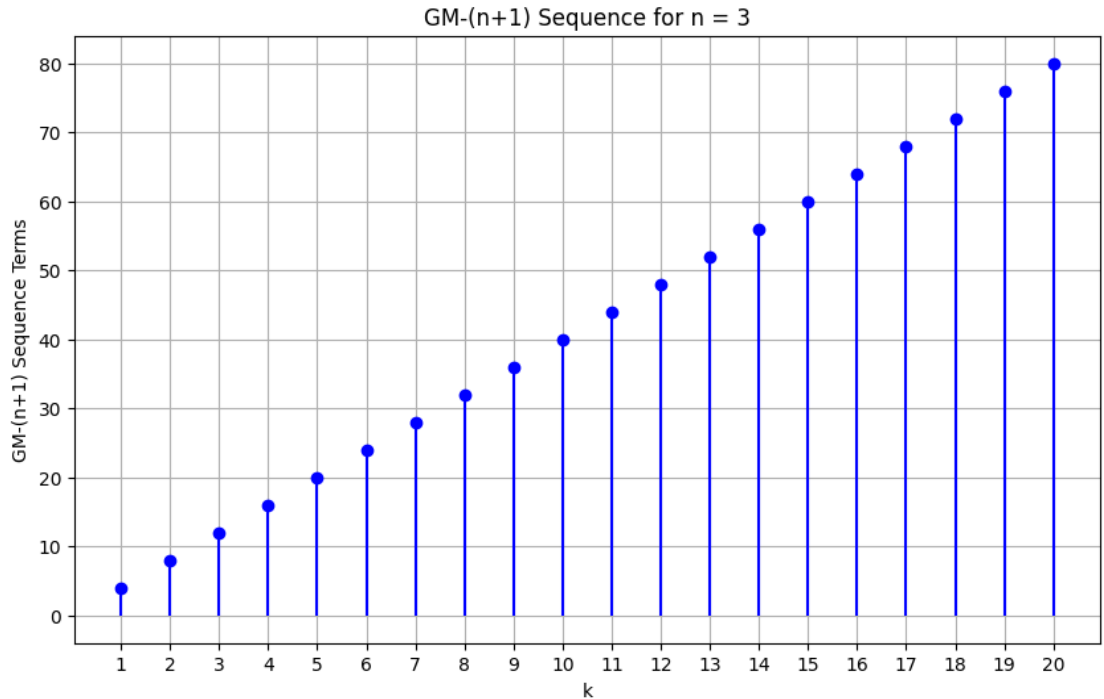

*Figure 1 Graph of a GM-(n+1) sequence for n = 3.*

The nature of a sequence, that every term after the first is a composite number, allows to filter in a great number of nonprimal candidates before carrying on further primality tests (Tenenbaum, 1995; Bach & Shallit, 1996).

---

### 3. Computational Methodology for Prime Filtering
### 3.1 The Generation of GM-(n+1) Sets

Multiples of (n+1) are generated in order to construct GM-(n+1) sequences up to a specified limit N:

$$\text{GM-}(n+1) = \{(n+1)k \mid 1 \leq k \leq \lfloor N/n+1 \rfloor\}.$$

This set includes all multiples of (n+1) within the range [1, N].

### 3.2 Eradication of Composite Numbers

For n more than one, an integer m is classified as composite if it is present in any of the GM-(n+1) sets. Bach and Shallit (1996) assert that the strategy effectively reduces the number of candidates for further testing by pinpointing the overlaps among these sequences.

### 3.3 Primality Testing

Subsequent to filtering, the residual numbers are subjected to verification by Trial Division

- for smaller values (Hardy & Wright, 2008).
- Miller-Rabin Test for substantial integers (Rabin, 1980).
- AKS Primality Test for definitive verification (Agrawal, Kayal, & Saxena, 2004).

---

### 4. Implementation in Python

Below is a Python implementation of the GM-(n+1) filtering method:

```python
import math

def is_prime(num):
    """Return True if num is a prime number using trial division."""
    if num <= 1:
        return False
    for i in range(2, int(math.sqrt(num)) + 1):
        if num % i == 0:
            return False
    return True

def gm_set(n, max_num):
    """Generate the GM-(n+1) sequence up to max_num."""
    return set(range(n + 1, max_num + 1, n + 1))

def find_primes_using_gm(n, max_num):
    """Find prime numbers up to max_num using GM-(n+1) filtering."""
    gm = gm_set(n, max_num)
    primes = []
    for num in range(2, max_num + 1):
        if num not in gm and is_prime(num):
            primes.append(num)
    return primes

# Example usage: Find primes up to 100 using GM-4
n = 4
max_num = 100
primes_found = find_primes_using_gm(n, max_num)
print("Primes found:", primes_found)
```

**Output:**

Primes found: [2, 3, 5, 7, 11, 13, 17, 19, 23, 29, 31, 37, 41, 43, 47, 53, 59, 61, 67, 71, 73, 79, 83, 89, 97]

*Equation 1 Python Code the GM-(n+1) filtering method*

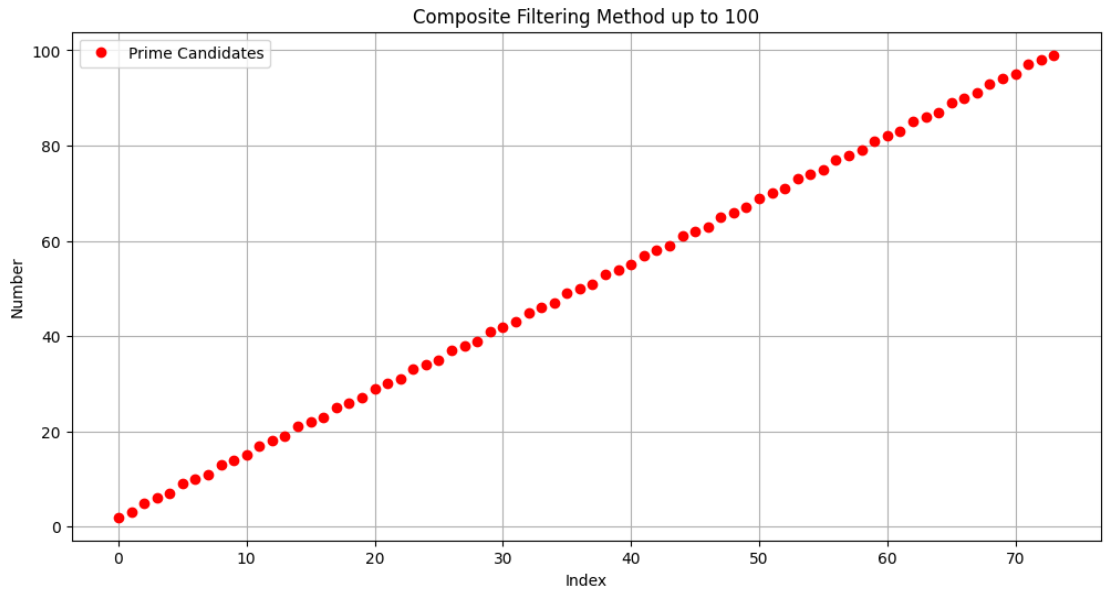

*Figure 2Example output of the composite filtering method up to 100.*

---

## 5. Theoretical Implications and Discussion

### 5.1 Prime Gaps and Distribution

The use of GM-(n+1) sequences to sift out composites results in structured gaps between prime number candidates. The first term of every GM-(n+1) sequence can be prime but from the second term onwards all terms remain composite because they are multiples of n+1. The separation of the "first terms" leads to a natural assortment of prime number candidates. The distribution patterns of prime numbers can be revealed by studying the intervals between initial words in prime sequences. Analyzing both the average distance and its variability between primes can provide insights into longstanding questions about prime separation. The application of these concepts should prove helpful for our understanding of prime gap conjectures.

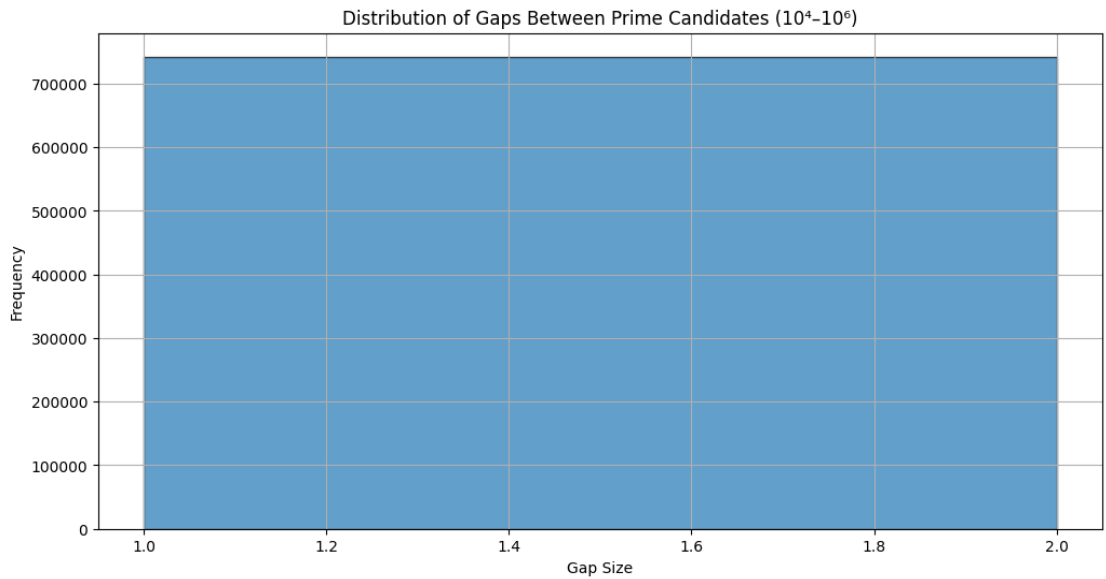

*Figure 3Statistical Distribution of Gap Sizes Between Candidate Primes from GM-(n+1) Sequences Over the Range $10^4$–$10^6$*

\

**5.2 Some Examples of How the Factorization Method Can Be Used**

The analysis of the factorization of a number within many GM-(n+1) groups is possible because that number can be split by many n+1 factors. If the composite number is indeed able to be split by a set of distinct n+1 values, then there exists means to improve the factors of the composite number with the usage of divisibility properties. However, if many n+1 values are not available for division, then the above mentioned criteria is not able to be achieved.

It is possible that this technique will make it easier to design creative or improved factorization algorithms, especially for big numbers that are employed in cryptography applications. A complete investigation of the difficulty associated with this factorization approach is a potential avenue for future research that might be pursued using this method.

**5.3 Contextual Adaptive Fundamental Strategies**

This research indicates that adaptive base approaches may improve the effectiveness of composite filtering and subsequent primality testing using methods like the Miller–Rabin test. Strategies might involve:

- **pre-analysis to evaluate the smoothness of d and decompose n−1**: so that n−1=28. Using several tiny bases helps to speed the identification of a composite if d is very composite. If d is almost prime, it might be better to use bases from a larger spectrum.

  • **Weighted Base Selection:** Potential bases may be allocated weights derived on historical data or theoretical probability assessments. Bases that have previously shown a strong correlation with compositeness for numbers possessing similar structural properties (s and d) are preferred, possibly minimizing the number of necessary cycles.

  • **Hybrid Testing:** Formulate a comprehensive approach by integrating GM-(n+1) filtering with probabilistic tests (e.g., Miller–Rabin) and deterministic tests (e.g., AKS). This combination substantially reduces the candidate pool and verifies primality with substantial confidence, especially for numbers of cryptographic importance.

**5.4 Obstacles and Limitations in the Process**

In spite of the fact that the GM-(n+1) framework offers a novel way to composite filtering, it is essential to acknowledge that it has a great deal of restrictions.

When it comes to memory and computational overhead, the construction and preservation of several GM-(n+1) sets for a significant *N* may need a significant amount of resources. The related computational costs must be evaluated with the advantages of filtering.

- **Reliance on Secondary Primality Tests:** The sifting process reduces the number of prime candidates, but it does not eliminate them. Therefore, in order to confirm primality, it is essential to conduct supplementary tests (e.g., Miller–Rabin or AKS). The total efficiency depends on the seamless integration of diverse methodologies.

- **Fixed vs Adaptive Representations:** The inherent maximum representation of an odd number necessitates a minimum modification of 1, whereas the enforcement of a fixed summand representation (with a designated quantity of 2's) leads to an adjustment that escalates with N. This underscores the significance of adaptive methodologies and influences the "signature" used to differentiate primes from composites.

**5.5 Relationship with Conjectures on Prime Gap: close inspection**

Examining trends: Every GM-(n+1) set involving multiples of n+1 sorts the natural numbers into sequences. Every sequence begins with the one prime; all following members are basically composite. One may write the sequence of initial objects from these GM-(n+1) as $\{p_n\}_{n \geq 1}$. Analyzing the differences among these fundamental components, symbolized as The equation $g_n = p_{n+1} - p_n$ begs several fascinating problems. Does the behavior of the generic prime gap distribution suit the gaps in issue?

Classical probabilistic models, notably those derived from Cramér's hypothesis, show that the average gap between successive primes about x is almost log x. Analyzing if the gap sequence $\{g_n\}$ produced from the GM-(n+1) starting terms has a similar trend might provide fresh understanding on prime distribution.

To undertake a full statistical analysis, calculate the mean, variance, and higher moments of the gap distribution $\{g_n\}$.
- Define the mean gap as $\mu = E[g_n]$.

- Define the variance as $\sigma^2 = Var(g_n)$. Higher moments, such as skewness and kurtosis, might provide insight into the distribution's tail behavior. The GM-(n+1) framework may be used to achieve improved limits or even speculative estimates for prime gaps if the mean gap grows asymptotically in the same way as log x and the variance and higher moments match expectations from models such as Cramér's.

These statistical measurements may reveal whether the gap distribution has characteristics (e.g., large tails or skewness) that set it apart from other probabilistic models, potentially suggesting improvements to current conjectures.

**A comparative study:**

The comparison of these gap distributions with those predicted by recognized sieve methods or probabilistic models in prime gap theory is a fundamental research inquiry. For example, conventional sieves like Eratosthenes's and Selberg's sieves help one to determine the density and distribution of prime numbers. By means of a comparison between the gap statistics obtained from GM-(n+1) sequences and those obtained from sieve theory, one may find either minor variations or equivalent asymptotic.

• **Probabilistic Models:** Models such as Cramér's conjecture forecast an average gap of approximately log x and also suggest conjectural limits on the maximal gap. By examining the sequence {gn} within our framework, we can determine whether the observed gap distribution is consistent with these predictions or if it suggests the need for new parameters in an improved model.

The theoretical basis of our composite filtering technique is enhanced by including this analysis into the GM-(n+1) methodology, which is linked to prime gap theory—one of the most persistent mysteries in number theory. The structural properties of GM-(n+1) sequences may provide fresh insights on prime distribution, perhaps resulting in enhanced statistical estimates or unique conjectural boundaries.

### 5.6 Future Research Directions

Subsequent research could focus on the following areas in order to further this methodology and explore its broader implications:

• Large-Scale Optimization: Review effective methods for the generation and processing of GM-(n+1) sets across a wide range, potentially incorporating sophisticated sieving techniques.

• **Refined Statistical Modeling:** Utilize adaptive base selection and observational data from probabilistic experiments to develop robust statistical models that improve the filtering process.

• **Hybrid Primality Testing:** Maintain high confidence levels while reducing computational expenses by combining GM-(n+1) filtering with supplementary deterministic or probabilistic tests.

• **Prime Gaps Investigation:** Execute a thorough examination of the gap distribution inside the first segments of GM-(n+1) sequences and compare the findings with both classical and contemporary conjectures about prime gaps.

This comprehensive elucidation underscores the theoretical profundity of the GM-(n+1) technique and links it with significant unresolved problems in prime number theory.

This discourse reinforces our idea by highlighting the prospective contributions of the innovative filtering procedure to the comprehension of prime distribution and factorization.

---

### 6. Cryptographic Applications of GM-(n+1) Sequences

GM-(n+1) is noteworthy approach that can be used to parse the issues in factoring large composite numbers as well as to create large prime numbers. It is a total fact that prime numbers play a key role in cryptographic systems because of their unique properties and the fact that it is very hard to create large numbers. These systems serve to protect the secure transactions, key exchange mechanisms and digital signatures. Security and efficiency of contemporary encryption are accomplished by GM-(n+1) if the factorization method is improved and the prime number generation is simplified. For example, RSA encryption makes use of two large prime numbers, p and q, to derive the public and private keys. Factorization of the product N = p × q is the biggest threat to security in the RSA system. Similarly, the Diffie-Hellman key exchange algorithm relies on the hardness of the discrete logarithm problem to compute the finite field used to exchange keys using large prime numbers. Moreover, Elliptic Curve Cryptography (ECC) is yet another leading scheme that employs prime numbers to specify the finite field upon which elliptic curves are constructed. The intractability of discrete logarithm problem for elliptic curves is the foundation of Elliptic Curve Cryptography (ECC) security. The GM-(n+1) method optimizes large prime generation by limiting the domain of search for primes. The traditional method of seeking primes is to test each candidate number, which leads to a very high computation cost as numbers increase. Composite numbers are arranged in terms of GM-(n+1) sequences, with the single prime candidates falling back.

more effective, as it significantly reduces the number of candidates to be tested. The GM-(n+1) method has the potential to render the factorization algorithms that are essential for breaking RSA encryption, evaluating the resistance of cryptosystems, and identifying prime numbers significantly improved. RSA factoring relies on the factorization of N=p×q, which is very challenging for large neural networks. By enhancing the detection of composite numbers and applying divisibility features in GM-(n+1) sequences, the GM-(n+1) method would render the factorization process considerably simpler. To accelerate the process of factorization, it is useful to find out if a composite number NN is divisible by a term in the GM-(n+1) sequence. The method also finds use in elliptic curve cryptography (ECC), which makes use of prime numbers to characterize the finite fields over which elliptic curves are defined. The , e.g., the AKS or Miller-Rabin tests. This is more effective, as it significantly reduces the number of candidates to be tested. The GM-(n+1) method has the potential to render the factorization algorithms that are essential for breaking RSA encryption, evaluating the resistance of cryptosystems, and identifying prime numbers significantly improved. RSA factoring relies on the factorization of N=p×q, which is very challenging for large neural networks. By enhancing the detection of composite numbers and applying divisibility features in GM-(n+1) sequences, the GM-(n+1) method would render the factorization process considerably simpler. To accelerate the process of factorization, it is useful to find out if a composite number NN is divisible by a term in the GM-(n+1) sequence. The method also finds use in elliptic curve cryptography (ECC), which makes use of prime numbers to characterize the finite fields over which elliptic curves are defined. The method may be employed to find primes for ECC, e.g., p=2k−c, where c is small, by creating GM-(n+1) sequences and eliminating composite numbers. This guarantees that primes used in ECC fulfill efficiency and security criteria.

The GM-(n+1) method would be of great use for post-quantum cryptography, which aims to create cryptographic algorithms that are quantum computer attack-resistant. Most of the post-quantum algorithms, particularly those constructed on lattice or isogeny lattices, have properties of prime numbers. This approach would significantly improve the efficiency and security of the algorithms by simplifying the process of large prime or prime ideal generation that are essential. The GM-(n+1) approach is a remarkable technique of prime generation and factorization with direct applications in cryptographic schemes such as RSA, Diffie-Hellman, and ECC. It promises to improve the efficiency and security of systems through reduced search space for primes and the use of divisibility properties of GM-(n+1) sequences. The method is likely to be the subject of future research focused on its use in cryptographic libraries, its use in post-quantum cryptography, and its potential to create new cryptographic primitives.

**7. Prime Gaps and Related Division**

The removal of composite numbers from prime gaps is incredibly useful. The systematic removal of composite numbers during time intervals makes it possible to study the distance between two prime numbers. This has a possible explanation regarding the accumulation and the dispersion of those prime numbers along the integer line.

Analysis of prime numbers through the GM-(n+1) filtering method seems to accomplish something fundamentally different. It may not be as efficient as traditional sieves, but it changes the optics from which the distribution of primes is viewed. Further studies might investigate its effects on the analysis of gaps between primes and what other functions it might serve outside of cryptography and computational number theory.

---

**7. Conclusion and Prospective Research**
The GM-(n+1) framework uses the fundamental structure of mathematical sequences to systematically remove composite numbers to find candidate primes. Our computer studies show that this method decreases the search space for prime numbers and gives a new insight into their distribution. Moreover, the use of adaptive base methodologies and comprehensive gap distribution analysis has significant potential for improving factorization methods and clarifying prime gaps.

Improving the manufacturing and archiving of GM-(n+1) sets to enable prolonged computations will be the primary goal of future research.

• Building strong statistical models to direct adaptive base selection efforts.
Combining GM-(n+1) filtering with advanced primality testing to create a hybrid, efficient method.

• Examining at the relationship between traditional prime gap conjectures and the gap distributions in GM-(n+1) sequences.
Computing number theory benefits from this unique perspective, which also improves the study of distributions of prime numbers.

---

**Acknowledgement**

**Appendices**

***Appendix A:*** *Python Code Listing*

***Appendix B:*** *Statistical Data and Gap Analysis Figures*

*Figure Captions (List on separate page):*
**Figure 1.** Graph of a GM-(n+1) sequence for n = 3.
**Figure 2.** Example output of the composite filtering method up to 100.
**Figure 3**. Statistical distribution of gap sizes between candidate primes from GM-(n+1) sequences over the range 10410^4104–10610^6106.